\documentclass{conm-p-l}
\usepackage{amssymb,latexsym, amscd}
\usepackage[all]{xy}

 \newtheorem{theorem}{Theorem}[subsection]

 \newtheorem{proposition}[theorem]{Proposition}
 \theoremstyle{definition}
 \newtheorem{definition}[theorem]{Definition}
 \theoremstyle{definition}
 
 \theoremstyle{remark}
 
 \numberwithin{equation}{subsection}

\usepackage{latexsym}
\usepackage{amssymb}

\newcommand{\ben}{\begin{equation}}
\newcommand{\een}{\end{equation}}

\newcommand{\integer}{\ensuremath{{\mathbb Z}}}

\newcommand{\real}{\ensuremath{{\mathbb R}}}
\newcommand{\complex}{\ensuremath{{\mathbb C}}}

\newcommand{\rational}{\ensuremath{{\mathbb Q}}}

\newcommand{\GL}[1]{\ensuremath{{\mathrm {GL}_{ #1 }}}}

\newcommand{\Aa}{{\mathcal A}}

\newcommand{\EE}{{\mathcal E}}
\newcommand{\FF}{{\mathcal F}}

\newcommand{\OO}{{\mathcal O}}
\newcommand{\CC}{\mathcal{C}}
\newcommand{\II}{\mathcal{I}}
\newcommand{\LL}{\mathcal{L}}
\newcommand{\MM}{\mathcal{M}}
\newcommand{\HH}{\mathcal{H}}

\newcommand{\Hom}{\mathrm{Hom}}
\newcommand{\Map}{\mathrm{Map}}
\newcommand{\End}{\mathrm{End}}
\newcommand{\Aut}{\mathrm{Aut}}

\newcommand{\To}{\longrightarrow}

\newcommand{\redu}{\mathrm{red}}
\newcommand{\Img}{\mathrm{Im\ }}
\newcommand{\Ker}{\mathrm{Ker\ }}
\newcommand{\per}{\mathrm{per}}
\newcommand{\even}{\mathrm{even}}
\newcommand{\odd}{\mathrm{odd}}
\newcommand{\rank}{\mathrm{rank}}
\newcommand{\alg}{\mathrm{alg}}
\newcommand{\topo}{\mathrm{top}}
\newcommand{\Res}{\mathrm{Res}}
\newcommand{\Spec}{\mathrm{Spec}}
\newcommand{\ch}{\mathrm{ch}}
\newcommand{\Amod}{{A-\mathbf{mod}}}
\newcommand{\op}{\mathrm{op}}

\newcommand{\Fr}{\mathrm{Fr}}
\newcommand{\md}{\,\,\mathrm{mod}\,\,}

\begin{document}

\title[XI Solomon Lefschetz Memorial Lecture Series]{XI Solomon Lefschetz Memorial Lecture Series:
Hodge structures in non-commutative geometry. \\ \smaller\smaller{(Notes by Ernesto Lupercio)} }
\author{Maxim Kontsevich}

\address{IHES, 35 route de Chartres, F-91440, France}

 \email{maxim@ihes.fr}

\begin{abstract}
Traditionally, Hodge structures are associated with complex projective varieties.
In my expository lectures I discussed a non-commutative generalization of Hodge structures
in deformation quantization and in derived algebraic geometry.
\end{abstract}

\maketitle

\section{Lecture 1. September 8th, 2005}

\subsection{} This talk deals with some relations between algebraic geometry and non-commutative
geometry, in particular we explore the generalization of Hodge structures to the non-commutative realm.

\subsection{Hodge Structures.}  Given
a smooth projective variety $X$ over $\complex$ we have a naturally defined pure Hodge structure (HS) on its cohomology, namely:
\begin{itemize}
\item $H^n(X,\complex) = \bigoplus_{p+q=n,\ p,q\geq 0} H^{p,q}(X)$ by considering $H^{p,q}$ to be the cohomology
represented by forms that locally can be written as $$\sum
a_{i_1,\dots,i_p;j_1,\dots,j_q} dz_{i_1} \wedge dz_{i_2} \wedge \ldots \wedge dz_{i_p}
\wedge d\bar{z}_{j_1} \wedge d\bar{z}_{j_2} \wedge \ldots \wedge d\bar{z}_{j_q}$$
\item $H^n(X,\complex)$ is the complexification  $H^n(X,\integer)\otimes\complex$ of a lattice of finite rank.
\item $H^{p,q} = \overline{H^{q,p}}$.
\end{itemize}

\subsection{} To have this Hodge structure (of weight $n$) is the same as having the decreasing
filtration $$F^p H^n := \bigoplus_{p' \geqslant p,\ p'+q'=n} H^{p',q'}$$ for we have $H^{p,q} = F^p H^n \cap \overline{F^q} H^n.$

\subsection{} What makes this Hodge structure nice is that whenever we have a family $X_t$ of varieties
algebraically dependent  on a parameter $t$ we obtain a bundle of
cohomologies with a flat (Gauss-Manin) connection $$H_t^n =
H^n(X_t,\complex)$$ over the space of parameters, and $F^p_t$ is a
{\it holomorphic} subbundle (even though $H^{p,q}_t$ is not).
Deligne developed a great theory of mixed Hodge structures that
generalizes this for any variety perhaps singular or non-compact.

\subsection{Non-commutative geometry.} Non-commutative geometry (NCG) has been developed by
Alain Connes \cite{Connes} with applications regarding foliations, fractals and quantum spaces in
mind, but not algebraic geometry. In fact it remains unknown what a good notion of non-commutative complex manifold is.

\subsection{} There is a calculus associated to NC spaces. Suppose $\bf k$ is a field and for the first talk the field will
always be $\complex$, only in the second talk finite fields become relevant. Let $A$ be a unital, associative algebra over $\bf k$. An idea
of Connes is to mimic topology, namely differential forms, and the de Rham differential in this framework. We define the Hochschild
complex $C_\bullet (A,A)$ of $A$ as a negatively graded complex (for we want to have all differentials of degree $+1$),
$$ \stackrel{\partial}{\longrightarrow} A\otimes A\otimes A \otimes A \stackrel{\partial}{\longrightarrow}
A \otimes A \otimes A \stackrel{\partial}{\longrightarrow} A\otimes A \stackrel{\partial}{\longrightarrow} A,$$
where $A^{\otimes k}$ lives on degree $-k+1$.
The differential $\partial$ is given by
$$\partial(a_0 \otimes \cdots \otimes a_n) = a_0 a_1 \otimes a_2 \otimes \cdots \otimes a_n - a_0 \otimes a_1 a_2 \otimes \cdots \otimes a_n $$
$$+ \ldots + (-1)^{n-1} a_0 \otimes a_1 \otimes \cdots \otimes a_{n-1} a_n + (-1)^n a_n a_0 \otimes a_1 \otimes \cdots \otimes a_{n-1}.$$
This formula is more natural when we write the terms cyclically:
\begin{equation}
\begin{array}{ccccccc} &  &  & a_0 &  &  &  \\
&  & \otimes &  & \otimes &  &  \\  & a_n &  &  &  & a_1 &  \\  & \otimes
 &  &  &  &  \otimes & \\  & \vdots &  &  & &   \vdots  & \\ &  & \otimes & & \otimes &  & \\  & & & a_i & & &\end{array}
\end{equation}
for $a_0 \otimes \cdots \otimes a_n$. It is very easy to verify that $\partial^2=0$.

\subsection{} The homology of the Hochschild complex has an abstract meaning
 $$\mathrm{Ker\ }\partial/\mathrm{Im\ }\partial = \mathrm{Tor}_\bullet^{A\otimes_k A^{\mathrm{op}} - \mathrm{mod}}(A,A).$$

\subsection{} An idea in NC geometry is that as $A$ replaces a commutative space the Hochschild homology of $A$ replaces in
turn the complex of differential forms.

\begin{theorem}[Hochschild-Konstant-Rosenberg, 1961, \cite{HKR}] Let $X$ be a smooth affine algebraic  variety, then if $A=\OO(X)$
we have $$HH_i(X):=H^{-i}(C_\bullet(A,A);\partial) \cong \Omega^i(X)$$ where $\Omega^i(X)$ is the space of $i$-forms on $X$.
\end{theorem}

The proof is very easy: consider the diagonal embedding $X\stackrel{\Delta}{\longrightarrow} X \times X$ and by
remembering that the normal bundle of $\Delta$ is the tangent bundle of $X$ we have
$$HH_\bullet(X) = \mathrm{Tor}_\bullet^{{\mathrm{Quasi-coherent}}(X\times X)} (\OO_\Delta,\OO_\Delta)$$ this together with a
local calculation gives the result.

\subsection{} The Hochschild-Konstant-Rosenberg theorem motivates us to think of $HH_i(A)$ as a space of differential forms of degree $i$
on a non-commutative space.

Note that if $A$ is non-commutative we have
$$H^0(C_\bullet(A,A);\partial)=A/[A,A].$$ Also, for commutative
$A=\OO(X)$,  given an element $a_0 \otimes\cdots \otimes a_n$ in
$C_\bullet (A,A)$ the corresponding form is given by $\frac{1}{n!}
a_0 da_1\wedge \ldots \wedge da_n.$

\subsection{} There is a reduced version of the complex $C^\redu_\bullet(A,A)$ with the same cohomology obtained by reducing modulo
constants all but the first
factor $$\longrightarrow A \otimes A/({\bf k}\cdot 1) \otimes A/({\bf k}\cdot 1)
\longrightarrow A \otimes A/({\bf k}\cdot 1) \longrightarrow A.$$

\subsection{} Connes' main observation is that we can write a formula for an additional differential
$B$ on $C_\bullet(A,A)$ of degree $-1$, inducing a differential on
$HH_\bullet(A)$ that generalizes the de Rham differential:
$$B(a_0\otimes a_1 \otimes \cdots \otimes a_n) = \sum_\sigma (-1)^\sigma 1 \otimes a_{\sigma(0)}
\otimes \cdots \otimes a_{\sigma(n)}$$
where $\sigma\in \integer/(n+1)\integer$ runs over all cyclic permutations. It is easy to verify
that $$B^2=0,\ \ \ B\partial +\partial B = 0, \ \ \  \partial^2=0,$$ which we depict pictorially as
$$\xymatrix{
                   \cdots \ar@/_/[rr]_\partial && A \otimes A/1 \otimes A/1 \ar@/_/[ll]_B \ar@/_/[rr]_\partial &&
                   A \otimes A/1 \ar@/_/[ll]_B \ar@/_/[rr]_\partial && A \ar@/_/[ll]_B
                  }$$
 and  taking cohomology gives us a complex $(\Ker \partial / \Img
\partial ; B)$. A naive definition on the de Rham cohomology in
this context is the homology of this complex $\Ker B / \Img B.$

\subsection{} We can do better by defining the negative cyclic complex $C^-_\bullet(A)$, which is formally a
projective limit (here $u$ is a formal variable, $\deg(u)=+2$):
$$C_\bullet ^-:= (C_\bullet^\redu(A,A)[[u]] ; \partial + u B) =
\lim_{\stackrel{\longleftarrow}{N}}(C_\bullet^\redu(A,A)[u]/u^N;\partial + u B).$$

\subsection{} We define the periodic complex as an inductive limit
$$C_\bullet ^\per:= (C_\bullet^\redu(A,A)((u)) ;
\partial + u B) = \lim_{\stackrel{\longrightarrow}{i}}(u^{-i} C_\bullet^\redu(A,A)[[u]];\partial + u B).$$
This turns out to be a ${\bf k}((u))$-module and this implies that
the multiplication by $u$ induces a sort of Bott periodicity. The
resulting cohomology groups called (even, odd) periodic cyclic
homology and are written (respectively) $$HP_\even(A), \ \ \ \ \
HP_\odd(A).$$ This is the desired replacement for de Rham
cohomology.

\subsection{} Let us consider some examples. When $A=C^\infty(X)$ is considered as a nuclear Fr\'echet algebra,
and if we interpret the symbol $\otimes$ as the topological tensor product then we have the canonical isomorphisms:
$$ HP_\even(A) \cong H^0(X,\complex) \oplus H^2(X,\complex) \oplus \cdots $$
$$ HP_\odd(A) \cong H^1(X,\complex) \oplus H^3(X,\complex) \oplus \cdots $$
\begin{theorem}[Feigin-Tsygan, \cite{Tsygan}] If $X$ is a affine algebraic variety (possibly \emph{singular})
and $X_\topo$ its underlying topological space then
$$HP_\even(A) \cong H^\even(X_\topo,\complex)$$ and $$HP_\odd(A) \cong H^\odd(X_\topo,\complex)$$
(these spaces
are finite-dimensional).
\end{theorem}
There is a natural lattice $H^\bullet(X,\integer)$ but we will see later that the ``correct'' lattice should be
slightly different.

\subsection{} Everything we said before can be defined for a differential graded algebra (dga) rather than only
for an algebra $A$. Recall that a dga $(A,d)$ consists of
\begin{itemize}
\item $A=\bigoplus_{n\in\integer} A^n$ a graded algebra with a graded product
$$A^{n_1} \otimes A^{n_2} \longrightarrow A^{n_1+n_2}.$$
\item $d_A \colon A^n \longrightarrow A^{n+1}$ a differential satisfying the graded Leibniz rule.
\end{itemize}
For example given a manifold $X$ on has  the de Rham dga $(\Omega^\bullet(X);d)$.

\subsection{} The definition of the degree for
$C_\bullet(A,A)$ is given by $$\deg(a_1\otimes\cdots\otimes a_n):= 1-n + \sum_i \deg(a_i).$$
It is not hard to see that $$\rank(HP_\bullet (A)) \leqslant \rank(HH_\bullet(A)),$$ and therefore,
if the rank of the Hochschild homology is finite so is the rank of the periodic cyclic homology.

\subsection{Hodge filtration on $HP_\bullet(A)$.} We define $F^n HP_\even(A)$ as the classes
represented by sequences $\gamma_i \in C_i(A,A)$, $i\in 2 \integer$ (namely $\sum \gamma_i u^{i/2}$) such
that $i \geq 2n$. Similarly we define $F^{n+1/2} HP_\odd(A)$.

\subsection{} We have an interesting instance of this situation in ordinary topology $A=C^\infty(X)$. Here we have:
$$HP_\even(A)=H^0(X)\oplus H^2(X)\oplus H^4(X)\oplus H^6(X) \oplus \cdots$$ and
$F^0 = HP_\even(A)$, $F^1=H^2(X)\oplus H^4(X)\oplus H^6(X) \oplus \cdots$ and so on.  In non-commutative
geometry this filtration
is the best you can do, for there is no individual cohomologies.

\subsection{} Let $X$ be an algebraic variety (not necessarily affine). Weibel \cite{Weibel97} gave
a sheaf-theoretic definition of $HH_\bullet(X)$ and
$HP_\bullet(X)$. Namely, if $X$ is covered by affine open charts
$$X=\bigcup_{1\leqslant i \leqslant r} U_i,$$ we obtain not an
algebra, but a cosimplicial algebra $$\Aa_k :=
\oplus_{(i_0,\dots,i_k)}\OO(U_{i_0}\cap\ldots\cap U_{i_k}),$$
whose total complex
$\mathrm{Tot}(C_\bullet(\Aa_\bullet))=C_\bullet(X)$ still has two
differentials $B$ and $\partial$ as before. In fact
$$H^\bullet(C_\bullet(X),\partial)=\mathrm{Tor}_\bullet^{{\mathrm{Quasi-coherent}}(X\times
X)} (\OO_\Delta,\OO_\Delta)$$
   is graded in both positive and negative degrees.

Weibel observed that one can recover a tilted version of the Hodge
diamond in this manner. For a smooth projective $X$ one has
 $$HP_\bullet(X):=HP_\bullet(\Aa_\bullet) = H^\bullet(X)$$ and the filtration we defined becomes
  $$F^i(HP_\bullet)=\bigoplus_{p=i+n/2} F^p H^n(X), \ \ \ i\in \frac{1}{2}\integer,$$ reshuffling thus
  the usual Hodge filtration.

Observe that in this example we have: $$ H^{p,q} \subset F^{\frac{p-q}{2}}. $$

\subsection{} In general one can directly replace an algebraic variety by a dga using a theorem by Bondal and Van den Bergh.
For example,    let $\EE$ be a sufficiently ``large'' bundle over
a smooth projective variety $X$.  You may take for instance $\EE =
\OO(0) + \OO(1) + \cdots + \OO(\dim X)$. Take the algebra $\Aa$ to be
$$\Aa := (\Gamma(X, \End(\EE) \otimes \Omega^{0,1});
\bar\partial).$$ Then one can show that one can repeat the
previous constructions obtaining the corresponding filtration.
Namely, the periodic cyclic homology (and the Hodge filtration) of
the dga $\Aa$ coincides with those of $X$.

\subsection{} Take $X_\alg$ to be a smooth algebraic variety over $\complex$ and let $X_{C^\infty}$ be
its underlying smooth manifold.
Consider the natural map $$X_{C^\infty} \longrightarrow X_\alg.$$ This map induces an isomorphism
$$ HP_\bullet(X_{C^\infty}) \longleftarrow
HP_\bullet(X_\alg).$$ This isomorphism is compatible with the Hodge filtrations but the filtrations are \emph{different}.

\subsection{} The next important ingredients are the integer lattices. Notice that $H^\bullet(X,\complex)$ has a natural
integer lattice $H^\bullet(X,\integer)$, which allows us to speak of periods, for example. There is also another lattice commensurable
 with $H^\bullet(X,\integer)$, namely, the topological $K$-theory $K^\bullet_{\topo}(X):=K^\even(X_{C^\infty})\oplus K^\odd(X_{C^\infty})$.

Let $A$ be an algebra, then we get $K_0(A)$ by considering the projective modules over $A$. There is a Chern character map
$$ \xymatrix{
        K_0(A) \ar[rr]^{\mathrm{ch}} \ar[rd]& & F^0(HP_\even(A)) \ar@{^{(}->}[r] & HP_\even(A) \\
& HC^-_0 (A) \ar[ur] & }$$ Here it may be appropriate to recall
that $HP_\bullet(A)$ is a Morita invariant and  therefore we can
replace
 $A$ by $A \otimes \mathrm{Mat}_{n\times n}$ for $\mathrm{Mat}_{n\times n}$ a matrix algebra.

If $\pi \in A$ is a projector (namely $\pi^2=\pi$) we have explicitly
$$\mathrm{ch}(\pi) = \pi - \frac{2!}{1!} (\pi-1/2)\otimes\pi\otimes\pi\cdot u + \frac{4!} {2!}
(\pi-1/2)\otimes\pi\otimes\pi\otimes\pi\otimes\pi\cdot u^2+\ldots.$$
There is a similar story for $K_1(A) \longrightarrow HP_\odd(A),$ and also for higher $K$-theory.

\subsection{} If $A=C^\infty (X)$ then the image of $K_0(A)=K^0_{\topo}(X)$ is up to torsion
$$H^0(X,\integer)\oplus H^2(X,\integer)\cdot 2\pi i \oplus H^4(X,\integer)\cdot (2\pi i)^2 \oplus \cdots.$$
We have of course $K^0_{\topo}(X)\otimes\rational =
H^\even(X,\rational)$ but the lattice is different and so Bott
periodicity is broken. In order to restore it we must rescale the
odd degree part of the lattice by the factor $\sqrt{2\pi i}$, and
then we obtain
$$H^1(X,\integer)\cdot\sqrt{2\pi i} \oplus H^3(X,\integer)\cdot(\sqrt{2\pi i})^3 \oplus \cdots$$

We call this new lattice the \emph{non-commutative integral cohomology}
$$H^\bullet_{\mathrm{NC}} (X, \integer) \subset HP_\bullet(C^\infty (X)).$$

\begin{proposition} For $A=C^\infty (X)$ the image up to torsion of
$$\mathrm{ch}:K_n(A) \longrightarrow HP_{(n\md 2)}(A)$$ is $$(2\pi i)^{n/2} H^{(n\md 2)}_{\mathrm{NC}} (X, \integer).$$
\end{proposition}

\subsection{} We are ready to formulate one of the \emph{main problems in non-commutative geometry.} Let $A$ be
a dga over $\complex$. The problem is to define a nuclear
Fr\'echet algebra $A_{C^\infty}$ satisfying Bott periodicity
$K_i(A_{C^\infty}) \cong K_{i+2}(A_{C^\infty}), i\ge 0$ together
with an algebra homomorphism $A \to A_{C^\infty}$ satisfying:
\begin{itemize}
\item The homomorphism $A \to A_{C^\infty}$ induces  an isomorphism $$HP_\bullet(A)\cong HP_\bullet(A_{C^\infty}),$$
\item $\mathrm{ch}:K_n(A_{C^\infty}) \longrightarrow HP_\bullet(A_{C^\infty})$ is a lattice, i.e. when we tensor
with $\complex$ we obtain an isomorphism
$$  K_n(A_{C^\infty})\otimes_\rational \complex \stackrel{\cong}{\longrightarrow} HP_{(n\md 2)}(A_{C^\infty}).$$
\end{itemize}

\subsection{} Consider for example the case of a commutative algebra $A$. Every such algebra is an inductive
limit of finitely generated algebras
 $$A = \lim_{\to} A_n$$ where each $A_n$ can be thought of as a singular affine variety. In general
 $HP_\bullet(A) \neq \lim_{\to} HP_\bullet(A_n)$, but the
 right-hand side is a better definition for $HP_\bullet(A)$. In this case we find that the lattice we are
 looking for is simply $$\lim_{\to}
 K^\bullet_\topo(\mathrm{Spec}\,A_n(\complex)).$$

\subsection{} We will attempt now to explain some non-commutative examples that are close to the commutative realm,
and are obtained by a procedure called \emph{deformation quantization.} Let us consider first a $C^\infty$ non-commutative space.

Let $T^2_\theta$ be the non-commutative torus (for $\theta \in \real$) so that $ C^\infty(T^2_\theta) $ is precisely all
the expressions of the form $$ \sum_{n,m\in\integer} a_{n,m} \hat{z}_1^n  \hat{z}_2^m, \ \ \ \ a_{n,m}\in\complex,$$ such
 that for all $k$ we have $$a_{m,n} = O((1+|n|+|m|)^{-k}),$$ and $$\hat{z}_1 \hat{z}_2 = e^{i \theta} \hat{z}_2 \hat{z}_1.$$

 For $\theta \in 2\pi \integer$ we get the usual commutative torus.

\subsection{} We will also consider some non-commutative \emph{algebraic} spaces obtained by deformation quantization.
Start by taking a smooth affine algebraic variety $X$ and a bi-vector field $\alpha \in \Gamma(X,\bigwedge^2 TX).$ Define,
as is usual, the
 bracket by $$\{f,g\}:= \langle \alpha , df\wedge dg\rangle.$$ The field $\alpha$ defines a Poisson structure iff the bracket satisfies
 the Jacobi identity. We will call $\alpha$ \emph{admissible at infinity} if there exists a smooth projective variety $\bar{X} \supset X$ and a
 divisor $\bar{X}_\infty = \bar{X} - X$ so that $\alpha$ extends to $\bar{X}$ and the ideal sheaf $\II_{\bar{X}_\infty}$
 is a Poisson ideal (closed under brackets).

\subsection{} The simplest instance of this is when $X=\complex^n$, $\bar{X}=\complex P^n$ and the admissibility
 condition for $\alpha = \sum_{i,j} \alpha_{i,j} \partial_i \wedge \partial_j$ reads $\deg(\alpha_{ij}) \leqslant 2$.

\subsection{} We have the following \cite{KontsevichDeformation}: \begin{theorem} If $X$ satisfies
 $$H^1(\bar{X},\OO) = H^2(\bar{X}, \OO) = 0,$$
 (e.g. $X$ is a rational variety) then there exists a canonical filtered algebra $A_\hbar$ over $\complex[[\hbar]]$
 (actually a free module over $\complex[[\hbar]]$) that gives a
 $*$-deformation quantization, and when we equal the deformation parameter to $0$ we get back $\OO(X)$. \end{theorem}

While the explicit formulas are very complicated the algebra $A_\hbar$ is \emph{completely canonical}.

\subsection{} This theorem raises the interesting issue of comparison of parameters. Take for example
the case $X=\complex^2$ and $\alpha=xy\frac{\partial}{\partial x}
\wedge \frac{\partial}{\partial y}.$ Here we can guess that
$$A_\hbar \cong \complex[[\hbar]] \langle \hat{X}, \hat{Y}\rangle
/ \left(\hat{X}\hat{Y}=e^\hbar \hat{Y} \hat{X}\right).$$ On the
other hand the explicit formula for $A_\hbar$ involves infinitely
many graphs and even for this simple example it is impossible to
get the explicit parameter $e^\hbar$. A priori one only knows that
just certain universal series $q(\hbar)=1+\dots$ should appear
with
 $\hat{X}\hat{Y}=q(\hbar) \hat{Y} \hat{X} $.

\subsection{} A slightly more elaborate example is furnished by considering (Sklyanin) \emph{elliptic algebras}. Here
we take $\bar{X}=\complex P^2$ and $\alpha=p(x,y) \frac{\partial}{\partial x} \wedge \frac{\partial}{\partial y}$
with $\deg(p)=3$. The divisor $\bar{X}_\infty \subset \complex P^2$ in this case is a cubic curve. We take
$X=\bar{X} - \bar{X}_\infty$, which is an affine algebraic surface and since it has a symplectic structure it also
has  a Poisson structure $\alpha$.

Its quantum algebra $A_\hbar$ depends on an elliptic curve $E$ and
a shift $x\mapsto x+x_0$ on $E$. The question is then: How to
relate $E$ and  $x_0$ to the bi-vector field $\alpha$ and the
parameter $\hbar$?

Again there is only one reasonable guess. Start with the bi-vector field $\alpha$ and obtain a $2$-form $\alpha^{-1}$
on $\complex P^2$ with a first order pole
at $\bar{X}_\infty$. Our guess is that $E=\bar{X}_\infty$. Taking residues we obtain a holomorphic
$1$-form $\Res(\alpha^{-1}) \in \Omega^1(E)$. The inverse of this $1$-form is a vector field $(\Res(\alpha^{-1}))^{-1}$
on $E$. Finally: $$x_0 = \exp \left( \frac{\hbar}{\Res(\alpha^{-1})}\right),$$ but to prove this directly seems to be quite challenging.

\subsection{} It is a remarkable fact that this comparison of parameters problem can be solved by considering the Hodge  structures.

\subsection{} Let us consider $X$ to be either a $C^\infty$ or an affine algebraic variety and $A$ to be $C^\infty(X)$
(respectively $\OO(X)$). The theory of deformation quantization
implies that all nearby non-commutative algebras and related
objects (such as $HP_\bullet$, $HH_\bullet$, etc.) can be computed
semi-classically. In particular nearby algebras are given by
Poisson bi-vector fields $\alpha$. Also $C_\bullet(A_\hbar,
A_\hbar)$ is quasi-isomorphic to the negatively graded complex
$(\Omega^{-i}(X), \LL_\alpha)$ where the differential is
$\LL_\alpha = [\iota_\alpha, d]$. If you want to see this over
$\complex[[\hbar]]$ simply consider the differential  $\LL_{\hbar
\alpha}$. We just described what Brylinski calls \emph{Poisson
homology}. The differential $B$ in this case is simply the usual
de Rham differential $B=d$. We would like to consider now
$HP_\bullet(A_\hbar)$. This is computed by the complex
$$(\bigoplus \Omega^i(X)[i][[\hbar]]((u)),ud+\hbar\LL_\alpha)\,,$$
 which is the sum of infinitely many copies of some finite-dimensional complex. Namely
 $HP_\bullet(A_\hbar)\hat{\otimes}_{\complex[[\hbar]]}\complex((\hbar))$ is $\complex((\hbar))$ tensored
 with the finite-dimensional cohomology of the $\integer/2$-graded complex:
$$\xymatrix{
                   \Omega^N  \ar@/_/[rr]_{\LL_\alpha} && \Omega^{N-1} \ar@/_/[ll]_d \ar@/_/[rr]_{\LL_\alpha} &&
                   \cdots  \ar@/_/[ll]_d \ar@/_/[rr]_{\LL_\alpha} && \Omega^{1}  \ar@/_/[ll]_d \ar@/_/[rr]_{\LL_\alpha}
                   && \Omega^0 \ar@/_/[ll]_d
                  }$$

\subsection{} We claim that the cohomology of this complex is $H^\bullet(X)$. The reason for this is really simple,
for we have that $$\exp(\iota_\alpha) d \exp(-\iota_\alpha) = d + [\iota_\alpha,d] + \frac{1}{2!}
[\iota_\alpha,[\iota_\alpha,d]]+\ldots = d+\LL_\alpha.$$
Here we used the fact that $[\alpha,\alpha]=0$ to conclude that only the first two terms survive.

\subsection{} Let us turn our attention to the lattice. Our definition uses $K_n(A)$ but we may get this lattice
by using the Gauss-Manin connection.
If we have a family of algebras $A_t$ depending on some parameter $t$, Ezra Getzler \cite{Getzler} defined a flat
connection on the bundle $HP_t$ over the parameter space. This allows us to start with the lattice
$\oplus_k H^k(X, (2\pi i)^{k/2} \cdot \integer)$ (up to torsion)
at $t:=\hbar=0$ in our situation. The parallel transport for the Gauss-Manin connection comes from the
above identification of periodic complexes given by the
 conjugation by $\exp (\hbar\iota_\alpha)$.

\subsection{} To compute the filtration we will assume that $X$ is symplectic, and therefore $\alpha$ is
non-degenerate. Again we set $\dim(X)=N=2n$, and
$\omega=\alpha^{-1}$ is a closed $2$-form. We also set $\hbar:=1$.
The following theorem is perhaps well known but in any case is
very simple:

\begin{theorem} For $(X,\omega)$ a symplectic manifold the Hodge filtration is given by
\begin{itemize}
\item $HP^\even$: $$F^{n-k/2} = e^\omega (H^0 \oplus \cdots \oplus H^{k}), \ \ k\in 2\integer$$
                              $$e^\omega H^0 \subset e^\omega(H^0\oplus H^2) \subset \cdots$$
\item $HP^\odd$: $$F^{n-k/2} = e^\omega (H^1 \oplus \cdots \oplus H^{k}), \ \ k\in 2\integer+1$$
                              $$e^\omega H^1 \subset e^\omega(H^1\oplus H^3) \subset \cdots$$
\end{itemize}
\end{theorem}

Notice that this is \emph{not} the usual Hodge filtration coming from topology:
$H^{2n}\subset H^{2n} \oplus H^{2n-2}\subset \cdots$

\begin{proof} Consider the $\integer/2$-graded complex
$$\xymatrix{
                   \Omega^{2n}  \ar@/_/[rr]_{\LL_\alpha} && \Omega^{2 n-1} \ar@/_/[ll]_d \ar@/_/[rr]_{\LL_\alpha} &&
                   \cdots  \ar@/_/[ll]_d \ar@/_/[rr]_{\LL_\alpha} && \Omega^{1}  \ar@/_/[ll]_d \ar@/_/[rr]_{\LL_\alpha}  &&
                   \Omega^0 \ar@/_/[ll]_d
                  }$$
where $\Omega^{k}$ lives in $F^{k/2}$.

The differential is not compatible with the filtration, nevertheless after the conjugation by
 $\exp(\iota_\alpha)$ we can use instead the complex $$\Omega^{2n} \stackrel{d}{\longleftarrow} \Omega^{2n-1}
 \stackrel{d}{\longleftarrow} \cdots \stackrel{d}{\longleftarrow} \Omega^{0}$$ and we would like to understand what
 happens to the filtration.

Let $*$ be the Hodge operator with respect to $\omega$ (the
Fourier transform in odd variables). Under this map the original
$\integer/2$-graded complex becomes
$$\xymatrix{
                   \Omega^{0}  \ar@/_/[rr]_{d} && \Omega^{ 1} \ar@/_/[ll]_{\LL_\alpha} \ar@/_/[rr]_{d} && \cdots
            \ar@/_/[ll]_{\LL_\alpha} \ar@/_/[rr]_{d} && \Omega^{2n-1}  \ar@/_/[ll]_{\LL_\alpha} \ar@/_/[rr]_{d}
            && \Omega^{2n} \ar@/_/[ll]_{\LL_\alpha}
                  }$$
where the filtration has been reversed.

The important remark here is that now $e^{\iota_\alpha}$ does preserve this filtration, transforming the last
complex into the complex
$$\xymatrix{
                   \ & \Omega^{0}  \ar[r]^{d} & \Omega^{1}  \ar[r]^{d} & \cdots \ar[r]^{d} & \Omega^{2n}\\
                   \deg & n & n-1/2 & \cdots & 0
                  }$$
Finally, all that remains to be seen is that $$e^{\omega\wedge\cdot} = e^{\iota_\alpha} * e^{-\iota_\alpha}.$$
\end{proof}

\subsection{} This theorem is related to the \emph{Lefschetz decomposition formula}\footnote{Lefschetz influence
in mathematics is clearly so large
 that is would be hard to give a talk in his honor without having the opportunity to mention his name at many points.}
 for K\"ahler manifolds. In fact we have obtained $$HP_\bullet(X_{\hbar \alpha}) = H^\bullet(X),$$ furthermore we have
 $$\lim_{\hbar\to 0} F^i HP_\bullet(X_{\hbar \alpha}) = F^i HP_\bullet(X)$$ if and only if we have a Lefschetz
 decomposition for the symplectic manifold.

If $(X,\omega)$ is K\"ahler compact we define the multiplication operator
 $$\omega\wedge : H^\bullet(X) \longrightarrow H^{\bullet+2}(X),$$ which is clearly nilpotent. The Lefschetz decomposition corresponds to
 the decomposition into Jordan blocks for this operator. In the case at hand the Lefschetz decomposition becomes the
 Hodge decomposition of a non-commutative space.

\subsection{} Consider $T^2_\theta$ the non-commutative torus. A result of Marc Rieffel \cite{Rieffel} states that
$T^2_\theta$ is Morita equivalent to $T^2_{\theta'}$
if and only if $$\theta' = \frac{a \theta +b}{c \theta + d},\ \ \ \
 \left(\begin{array}{cc}a & b \\c & d\end{array}\right)\in \mathrm{SL}_2(\integer).$$
If you consider in this case $$ HP_\even(T^2_\theta)=H^0\oplus H^2\longleftarrow K_0(T^2_\theta) $$
you can see that $K_0(T^2_\theta)$ contains the semigroup
of bona fide projective modules, producing a half-plane in the lattice bounded by a line of slope $\theta/(2\pi)$,
which from our point of view can be identified with the
Hodge filtration $F^1$. This helps to clarify the meaning of Rieffel's theorem. In this example we get an
interesting filtration only for $HP^\even$, and nothing for $HP^\odd$.

\subsection{} Consider an elliptic curve: $$E=\complex/\left(\integer+ \tau \integer\right),\ \ \ \  \Im(\tau)>0.$$
Here $$HP_\even(E) = H^0(T^2)\oplus H^2(T^2),$$ $$HP_\odd(E)=
H^1(T^2,\complex) \supset
 H^{1,0}(E)\,\,\{ \mathrm{ terms\ in\ the\ Hodge\ filtration}\}.$$
This becomes in terms of generators
$$\complex\otimes\left(\integer e_0 \oplus \integer e_1\right)
\supset \complex \cdot (e_0 + \tau e_1).$$ While in the
corresponding situation for $T^2_\theta$ the filtration can be
written as
$$\complex\otimes\left(\integer \tilde{e}_0 + \integer \tilde{e}_1\right) \supset
\complex \cdot (\tilde{e}_0 + \frac{\theta}{2\pi} \tilde{e}_1).$$
All this confirms a general belief that non-commutative tori are
limits of elliptic curves as $\tau\to \real$. Also $E$ can be seen
as a quotient of a 1-dimensional complex torus $\complex^\times$,
while $T^2_\theta$ plays the role of a \emph{real} circle modulo a
$\theta$-rotation.

\subsection{} I shall finish this lecture with one final puzzle. Namely, in the previous example the
Hodge filtrations do not fit. In the elliptic curve the
interesting Hodge structure detects parameters in odd cohomology while in the non-commutative torus
the Hodge filtration detects parameters in even cohomology.

A reasonable guess for the solution of this puzzle is that one should tensor by a simple super-algebra
(discovered by Kapustin in the Landau-Ginzburg model) given by
$$A=\complex[\xi]/(\xi^2=1)$$ with $\xi$ odd. Here $$HP_\bullet(A) = \complex^{0|1}.$$ The question
is: How does this super-algebra naturally arise from the limiting process $\tau\to\real$ sending an
elliptic curve to a foliation?

\section{Lecture 2. September 9th, 2005.}

\subsection{Basic Derived Algebraic Geometry.} This field started by A. Bondal and M. Kapranov in
Moscow around 1990. Derived algebraic geometry is much simpler
than algebraic geometry. While algebraic geometry starts with
commutative rings and builds up spectra via the Zariski topology
and the theory of sheaves, in derived algebraic geometry there is
no room for many of these concepts and the whole theory becomes
simpler.

\subsection{} Let us start by commenting on the algebraization of the notion of space. If we begin with a
(topological) space $X$, first one can produce an algebra $A=\OO(X)$, its algebra of functions. Next we
assign an abelian category to this algebra, the abelian category $\Amod$ of $A$-modules.

At every stage we \emph{insist} in thinking of the space as the remaining object: The category $\Amod$ \emph{is}
the space. The abelian category $\Amod$ has a nice subfamily, that of \emph{vector bundles}, namely finitely
generated projective modules. Recall that projective modules are images of $(n\times n)$-matrices $\pi : A^n \to A^n$ satisfying $\pi^2=\pi$.

The final step consists in producing from the category $\Amod$ a triangulated category $D(\Amod)$ that goes by
the name of the \emph{derived category}. $$X\mapsto \OO(X)=A \mapsto \Amod \mapsto D(\Amod).$$

\subsection{} While the abelian category $\Amod$ is nice we are still forced to keep track of whether a functor
is left-exact or right-exact, etc. This is greatly simplified in the derived category $D(\Amod)$.

The derived category $D(\Amod)$ is built upon infinite $\integer$-graded complexes of free $A$-modules and considering homotopies.

\subsection{} We take one step further  and consider the subcategory $\CC_X \subset D(\Amod)$ of \emph{perfect complexes.}
 A perfect complex is a finite length complex of finitely generated projective $A$-modules (vector bundles).
 All of the
 above can be generalized to dg algebras.

\subsection{} We are ready to make  important definitions:
\begin{definition}
Let $k$ be a field. A $k$-linear space $X$ is a small triangulated category $\CC_X$ that is Karoubi closed
(namely all projectors split), enriched by complexes of $k$-vector spaces.
 In particular for any two objects $\EE$ and $\FF$ we are given a complex $\Hom_{\CC_X}^\bullet(\EE,\FF)$ such that
  $$\Hom_{\CC_X}(\EE,\FF) = H^0(\Hom_{\CC_X}^\bullet(\EE,\FF)).$$
\end{definition}

\begin{definition} A $k$-linear space $X$ is algebraic if $\CC_X$ has a generator (with respect to taking cones and
direct summands). \end{definition}

\subsection{} The following holds:

\begin{proposition} The category $\CC_X$ has a generator if and only if there exists a dga $A$ over $k$ such
that $$\CC_X \cong \mathrm{Perfect}(\Amod).$$ \end{proposition}

This proposition allows us to forget about categories and consider
simply  dga-s (modulo a reasonable definition of derived Morita
equivalence).

\subsection{} There is a nice relation with the notion of scheme:

\begin{theorem}[Bondal, Van den Bergh \cite{Bergh}] Let $X$ be a scheme of finite type over $k$, then $\CC_X$
has a generator. \end{theorem}

The moral of the story in derived algebraic geometry is that \emph{all spaces are affine}.

\subsection{} The following example is due to Beilinson \cite{Beilinson}. Consider $X=\complex P^n$. Then
$$D^b(\mathrm{Coherent}(X))=\mathrm{Perfect}(X) = \mathrm{Perfect}(\Amod),$$ where $A=\End(\OO(0)\oplus\cdots\oplus\OO(n)).$
A finite complex of finite-dimensional representations of $A$ is the same as a finite complex of vector bundles over $X=\complex P^n$.

\subsection{} We make a few more definitions.

\begin{definition} An algebraic $k$-linear space $X$ is compact
 if for every pair of objects $\EE$ and $\FF$ in $\CC_X$ we have that $$\sum_{i \in \integer}\rank\ \Hom(\EE,\FF[i])  < \infty.$$
 In the language of the dga $(A,d_A)$ this is equivalent to: $$\sum_{i\in\integer} \rank\ H^i(A,d_A) <\infty.$$
\end{definition}

\begin{definition} We say that an algebraic $k$-linear space $X$ is smooth if
$$A\in \mathrm{Perfect}(A\otimes A^{\mathrm{op}} - \mathbf{mod}).$$
\end{definition}

\begin{definition} (a version of Bondal-Kapranov's) $X$ is saturated if it is smooth and compact.
\end{definition}

This is a good replacement for the notion of smooth projective variety.

\subsection{} The following concerns the moduli of saturated spaces:

\begin{proposition}[Finiteness Property] The moduli space of all saturated $k$-linear spaces $X$ modulo isomorphisms
can be written as a countable disjoint union of schemes of finite type:
$$ \coprod_{i\in I} S_i / \sim $$ modulo an algebraic equivalence relation.
\end{proposition}

\subsection{Operations with saturated spaces.} We have several basic operations inherited from the operations on algebras:
\begin{itemize}
\item[(i)] Given a space $X$ we can produce  its opposite space
$X^{\mathrm{op}}$ by sending the algebra $A$ to its opposite
$A^{\mathrm{op}}$. \item[(ii)] Given two spaces $X$ and $Y$ we can
define their tensor product $X\otimes Y$ by multiplication of
their corresponding dga's $A_X \otimes A_Y.$ \item[(iii)] Given
$X$, $Y$ we define the category $\Map(X,Y) := A_X^{\mathrm{op}}
\otimes A_Y - \mathbf{mod}.$ \item[(iv)] There is a nice notion of
gluing which is absent in algebraic geometry. Given $f: X\to Y$
(namely a $A_Y-A_X$-bimodule $M_f$) construct a new algebra $A_{X
\cup_f Y}$ by considering upper triangular matrices of the form
$$ \left(\begin{array}{cc} a_x & m_f \\0 & a_y \end{array}\right),$$ with $a_x \in A_X$, $a_y \in A_Y$ and $m_f\in M_f.$
\end{itemize}

Beilinson's theorem can be interpreted as stating that $\complex
P^n$ is obtained by gluing $n+1$ points. This does not sound very
geometric at first. An interesting outcome is that we get an
unexpected action of the braid group in such decompositions of
$\complex P^n$.

Notice that the cohomology of the gluing is the direct sum of the cohomologies of the building blocks.

\subsection{Duality theory} The story is again exceedingly simple for saturated spaces. There is a canonical Serre
functor $$S_X \in \Map(X,X)$$ satisfying $$\Hom(\EE,\FF)^* = \Hom(\FF,S_X(\EE))^*$$ In terms of the dga $A_X$ we have
that $$S_X^{-1}  = R\Hom_{A\otimes A^{\mathrm{op}}}(A,A\otimes A^{\mathrm{op}}).$$

In the commutative case this reads $$S_X=K_X[\dim X]\otimes\cdot\,\,,$$ where $K_X=\Omega^{\dim X}.$

\subsection{} It seems to be the case that there is a basic family of objects from which (almost)
everything can be glued up: Calabi-Yau spaces.
\begin{definition}
A \emph{Calabi-Yau saturated space} of dimension $N$ is a saturated space $X$ where the Serre functor $S_X$ is the shifting functor $[N]$.
\end{definition}

There are reflections of various concepts in commutative algebraic
geometry such as positivity and negativity of curvature in the
context of separated spaces. Here we should warn the reader that
sometimes it is impossible to reconstruct the commutative manifold
from its saturated space: several manifolds produce the same
saturated space. Think for example about the Fourier-Mukai transform.

\subsection{} We also have a $\integer/2$-graded version of this theory. We require all complexes and shift functors to be 2-periodic.

\subsection{Examples of saturated spaces.}
\begin{itemize}

\item Smooth proper schemes. They form a natural family of
saturated spaces.

\item Deligne-Mumford stacks that look locally like  a scheme $X$
with a finite group $\Gamma$ acting $X$. By considering (locally)
the algebra $A=\OO(X) \rtimes k[\Gamma]$ we can see immediately
that they also furnish examples of saturated spaces.

\item Quantum projective varieties. Suppose we start by
considering an ample line bundle $L$ over a smooth projective
variety $X$. Say we have $\alpha_X$ a bi-vector field defined over
$L-\mathbf{0}$ the complement of the zero section. We assume that
$\alpha_X$ is invariant under $\mathbb{G}_m = \GL{1}$. Deformation
quantization implies that we obtain a saturated non-commutative
space over $k((\hbar))$.

\item Landau-Ginzburg models. This is a $\integer/2$-graded example. Here we are given a map
$$f : X\To \mathbb{A}^1, \ \ \ f\in\OO(X), \ \ \ f\not\equiv 0,$$ where $X$ is a smooth non-compact variety
and $\mathbb{A}^1$ is the affine line. The idea comes from the
B-model in string theory. The category $\CC_{(X,f)}$ consists of
\emph{matrix factorizations}. In the affine case the objects are
super-vector bundles $\EE = \EE^\even \oplus \EE^\odd$ over $X$,
together with a differential $d_\EE$ such that $$d_\EE^2= f\cdot
\mathrm{Id}.$$ In local coordinates we are looking for a pair of
matrices $(A_{ij})$ and $(B_{ij})$ so that $$A\cdot B = f \cdot
\mathrm{Id}.$$

We define $\Hom((\EE,d_\EE),(\tilde{\EE},d_{\tilde{\EE}}))$ to be the complex of $\integer/2$-graded spaces
$\Hom_{\OO(X)}(\EE,\tilde{\EE})$ with differential
$$d(\phi) = \phi \cdot d_\EE - d_{\tilde{\EE}} \cdot \phi.$$ It is very easy to verify that $d^2=0.$

The generalization of the definition of $\CC_{(X,f)}$ to the global case is
 due to Orlov \cite{Orlov}. Consider $Z=f^{-1}(0)$ as a (possibly nilpotent) subscheme of $X$. Then
 $$\CC_{(X,f)} := D^b(\mathrm{Coherent}(Z))/\mathrm{Perfect}(Z).$$
 We expect $\CC_{(X,f)}$ to be saturated if and only if $X_0=\mathrm{Critical}(f) \cap f^{-1}(0)$ is compact.
 This is undoubtedly an important new class of triangulated categories.

\item A beautiful final example is obtained by starting with a
$C^\infty$ compact symplectic manifold $(X,\omega)$, with very
large symplectic form $[\omega]\gg 0$ (This can be arranged by
replacing $\omega$ by $\lambda\omega$ with $\lambda$ big.) The
\emph{Fukaya category} $\FF(X,\omega)$ is defined by taking as its
objects Lagrangian submanifolds and as its arrows holomorphic
disks with Lagrangian boundary conditions (but the precise
definition is not so simple.)

Paul Seidel \cite{Seidel} has proposed an argument showing that in
many circumstances $\FF(X,\omega)$ is saturated. This is a
manifestation of \emph{Mirror Symmetry} that says that
$\FF(X,\omega)$ is equivalent to $\CC_{(X,\omega)^\vee}$, where
${(X,\omega)^\vee}$ is the mirror dual to $(X,\omega)$. While
originally mirror symmetry was defined only in the Calabi-Yau
case, now we expect that the mirror dual to a general symplectic
manifold will be dual to a category of Landau-Ginzburg type. In
any case in many known examples the category is glued out of
Calabi-Yau pieces.

\end{itemize}

\subsection{} In derived algebraic geometry  there are a bit more spaces but much more identifications and symmetries
than in ordinary algebraic geometry. For instance, when $X$ is
Calabi-Yau then $\Aut(\CC_X)$ is huge and certainly much bigger
than $\Aut(X)$. Another example: there are two different $K3$
surfaces $X$, $X'$ that have the same $\CC_X=\CC_{X'}$, but in
fact $X$ and $X'$ need not be diffeomorphic. Again think of the
symmetries furnished by the Fourier-Mukai transform.

\subsection{Cohomology} Let us return to the subject of cohomology. First, let us make an important
 remark. If $A$ is a saturated dga then its Hochschild homology $H_\bullet(A,A):=H^\bullet(C_\bullet(A,A))$ is of
 \emph{finite rank}, and the rank
 of the periodic cyclic homology is bounded by $$\rank HP_\bullet(A) \leqslant \rank H_\bullet(A,A).$$ In the case
 in which $A$ is a
 commutative space we have $HP_\bullet(A)\cong H^\bullet_{\mathrm{de Rham}}(X)$ and
 $H_\bullet(A,A) \cong H^\bullet_{\mathrm{Hodge}}(X)$.

\subsection{} This motivates the following definition:
\begin{definition} For a saturated space $X$ over $k$ the Hodge to de Rham spectral sequence  is said to collapse if
$$ \rank\ H_\bullet(A,A) = \rank\ HP_\bullet(A) .$$
This happens if and only if for all $N \geqslant 1$ we have that $H^\bullet(C^\redu_\bullet(A,A)[u]/u^N,\partial+u B)$
is a free flat $k[u]/u^N$-module.
\end{definition}

\subsection{The Degeneration Conjecture:} For any saturated $X$ the Hodge to de Rham spectral sequence
collapses.\footnote{ See the very promising work of Kaledin that
has appeared since,
 \cite{Kaledin1, Kaledin2}.}

This conjecture is true for commutative spaces, for quantum projective schemes and for Landau-Ginzburg models $(X,f)$.

There are two types of proofs in the commutative case. The first
uses K\"ahler geometry and resolution of singularities. This
method seems very hard to generalize. The second method of proof
uses finite characteristic and Frobenius homomorphisms, this is
known as the \emph{Deligne-Illusie method} and we expect it (after
D.~Kaledin) to work in general.

\subsection{} Let us assume this conjecture from now on. We have then a vector bundle $\HH$ over $\Spec\ k[[u]]$ and we
will call $\HH_u$ the fiber. The space of sections of this bundle
is $HC^-_\bullet(A)$. This bundle carries a canonical connection
$\nabla$ with a first or second order pole at $u=0$.

In the $\integer$-graded case we have a $\mathbb{G}_m$-action, $$\lambda \in k^\times,\ \ \ \ u\mapsto \lambda^2 u,$$
defining the connection. The monodromy of the connection is $1$ on $HP^\even(A)$ and $-1$ on $HP^\odd(A)$. In this
case the connection has a first order pole at $u=0$ and the spectrum of the residue of the connection is $\frac{1}{2}\integer$.

The $\integer/2$-graded case is even nicer, for the connection can be written in a universal way with an explicit but
complicated formula containing the sum of five terms
 (see \cite{KontsevichSoibelman}). There is a reason for this connection to exist, and we explain it in two steps.

\begin{itemize}
\item Recall that if you have a family of algebras $A_t$ over a parameter space you get a flat connection on the
bundle $HP_\bullet(A_t)$ whose
formula tends to be very complicated.

\item Consider the moduli stack of $\integer/2$-graded spaces. We have an action of $\mathbb{G}_m$:
$$(A,d_A) \mapsto (A,\lambda d_A).$$ This corresponds in string theory to the renormalization group flow.
The fixed points of this action contain $\integer$-graded spaces (but there are also the elements of fractional charge,
and the quasihomogeneous singularities.)  This corresponds to a scaling $u \mapsto \lambda u$ and therefore produces the
desired connection. In this case we have a second order pole at $u=0$.
\end{itemize}

\subsection{} The basic idea is that the connection $\nabla$ replaces the Hodge filtration.

Notice that a vector space together with a $\integer$-filtration is the same as a vector bundle over $k[[u]]$ together with
a connection with first order pole at $u=0$ and with trivial monodromy. Of course, our connection is more complicated but it
is generalizing the notion of filtration.

\subsection{} We will use now the Chern character $$\ch : K_0(X) \To \{ \mathrm{covariantly\ constant\ sections\ of\ the\ bundle\ }\HH\}.$$
If in particular we consider the Chern class of $\mathrm{id}_X \in
\CC_{X \times X^\op}$ we have that $\ch(\mathrm{id}_X)$ is a
covariantly constant pairing $$\HH_u\stackrel{\cong}{\To}
\HH_{-u}^*$$ that is non-degenerate at $u=0$.

\subsection{} Let us now describe the construction of an algebraic model for a string theory of type IIB.

Let $X$ be a saturated algebraic space together with:
\begin{itemize}

\item A Calabi-Yau structure (this exists if the Serre functor is
isomorphic to a shift functor)\footnote{In a sense a Calabi-Yau
structure \emph{is}  more or less a choice of isomorphism between
Serre's functor and a shift functor.} To be precise a Calabi-Yau
structure is a section $\Omega_u \in HC^-_\bullet(A)=\Gamma(\HH)$
such that as $\Omega_{u=0}$ is an element in Hochschild homology
of $X$, that in turn gives a functional on $H_\bullet (A, d_A)$
making it into a Frobenius algebra, see
\cite{KontsevichSoibelman}.

\item A trivialization  of $\HH$ compatible with the pairing
$$\HH_u\otimes \HH_{-u} \To k$$ and so the pairing becomes constant.

\end{itemize}

If the main conjecture is true and the Hodge to de Rham
degeneration holds then from such $X$ we can construct a
2-dimensional cohomological quantum field theory. The state space
of the theory will be $\HH_0=HH_\bullet(X)$. As a part of the
structure one gets maps
$$\HH_0^{\otimes n} \To H^\bullet(\overline{\MM}_{g,n},k).$$ We will
not describe the whole (purely algebraic) construction here, but we shall just say that it provides solutions to holomorphic anomaly equations.

It seems to be the case that when we apply this procedure to the Fukaya category we recover the usual Gromov-Witten
invariants for a symplectic manifold. It is very interesting
to point out that the passage to stable curves is dictated by both, the degeneration of the spectral sequence, and the
trivialization of the bundle. This fact
 was my main motivation for the Degeneration Conjecture.

In the $\integer$-graded case a Calabi-Yau structure requires a
volume element $\Omega$ and a splitting of the non-commutative
Hodge filtration compatible with the Poincar\'e pairing.

\subsection{} We can make an important definition.
\begin{definition} A non-commutative \emph{ Hodge Structure} over $\complex$ is a holomorphic super vector bundle $\HH^{\mathrm{an}}$
over $D=\{|u|\leqslant 1,\ u\in\complex\}$ with connection
$\nabla$ outside of $u=0$, with a second order pole and a regular
singularity\footnote{By a regular singularity we mean that
covariantly constant sections grow only polynomially. Therefore
under a
 meromorphic gauge transformation we end up with a first order pole.}  together with a covariantly constant finitely generated
 $\integer/2$-graded abelian group $K_u^\topo$ for $u\neq 0 $ such that $K_u^\topo\otimes\complex = \HH_u^{\mathrm{an}}$.
\end{definition}

In the $\integer$-graded case the pole has order one, and the
lattice $K_u^\topo$ comes from the topological K-theory.

\subsection{The Non-commutative Hodge Conjecture.} Let $X$ be a saturated space. Consider the map
 $$\rational \otimes \mathrm{Image}(K_0(\CC_X) \To \Gamma(\HH^{\mathrm{an}}(X)))
  \to \rational \otimes \Hom_{\mathrm{NC-Hodge-structures} }
 (\mathbf{1},\HH^{\mathrm{an}}(X)).$$

The conjecture says that this map is an isomorphism.

One can introduce a notion of a \emph{polarized} NC Hodge structure.
The existence of a polarization in addition to the Hodge conjecture imply that the image of $K_0(\CC_X)$ is
$$K_0(\CC_X)/\mathrm{numerical\ equivalence},$$ where numerical equivalence is the kernel of a pairing
$$\langle,\rangle:K_0\otimes K_0 \To \real,$$ given by $$\langle \EE, \FF \rangle  = \chi (R\Hom(\EE,\FF)).$$
This pairing is neither symmetric nor anti-symmetric, so \emph{a priori} it could have left and right kernels, but the
Serre functor ensures us that they coincide.

\subsection{} We can now go on \`a la Grothendieck and define a category of non-commutative pure motives. Consider $X$
a saturated space over $k$. We define now
$$\underline{\Hom}(X,Y) = \rational \otimes K_0(\Map(X,Y))/\mathrm{numerical\ equivalence}.$$
Ordinarily one takes algebraic cycles of all possible dimensions
on the product of two varieties. In our situation we must be
careful to add direct summands. This should be equivalent to a
category of representations of the projective limit of some
reductive algebraic groups over $k$. This  non-commutative motivic
Galois group responsible for such representations  is much larger
than usual because of the $\integer/2$-gradings.

\subsection{} We can also discuss \emph{mixed motives} in this context. They are a replacement of Voevodsky's triangulated
category of mixed motives.
We start again by considering saturated spaces $X$ but now we want to define a new $\Hom(X,Y)$ space as the $K$-theory
spectrum of $\Map(X,Y)$ (an infinite loop space.)
We can canonically form the triangulated envelope. Notice that this construction contains ordinary mixed motives for usual
varieties modulo the tensoring by $\integer(1)[2]$.

\subsection{Crystalline cohomology and Euler functors.} Consider an algebra $A$ flat over $\integer_p$ (the $p$-adic integers)
and saturated over $\integer_p$. We expect a canonical Frobenius isomorphism
$$\Fr_p:H^\bullet (C_\bullet^\redu(A,A)((u)),\partial+uB) \stackrel{\cong}{\To}
H^\bullet(C^\redu_\bullet (A,A)((u)), \partial+puB),$$ as
$\integer_p((u))$-modules preserving connections.

Such isomorphism does exist in the commutative case, given a
smooth $X$ over $\integer_p$ we have
$$H^\bullet(X,\Omega_X;d) \cong H^\bullet(X,\Omega_X;pd)$$
for $p>\dim X$.

Using the holonomy of the connection $\nabla_{\frac{\partial}{\partial u}}$ we can go from $u$ to $pu$ and get an
operator $\Fr_p$ with coefficients in $\rational_p$.

We can state

\subsection{The non-commutative Weil conjecture.} Let $\lambda_\alpha \in \Spec Fr_p$ then
\begin{itemize}
\item $\lambda_\alpha \in
\overline{\rational}\subset\overline{\rational}_p.$ \item For all $\ell \neq
p$, then $|\lambda_\alpha|_\ell = 1$. \item
$|\lambda_\alpha|_\complex = 1$.
\end{itemize}

\subsection{} For the Landau-Ginzburg model $X=\mathbb{A}^1$ and $f=x^2$ we get that
the cohomology is 1-dimensional, $\lambda \in \rational_p$ and
$$\lambda = \left( \frac{p-1}{2} \right) !\,\,(\md p) ,$$
$$
\lambda^4=1.$$

The period for the Hodge structure is $\sqrt{2 \pi}$. There is a
reasonable hope for the existence of the Frobenius isomorphism
(cf. the work of D.~Kaledin.)

\subsection{} For every associative algebra $A$ over $\integer/p$ there is a canonical linear map $H_0(A,A) \To H_0(A,A)$
given by $$a \mapsto a^p.$$ Recall that $H_0(A,A)=A/[A,A].$ There
are two things to verify, that this map is well defined and that
it is linear. (This is a very pleasant exercise.) Moreover this
map lifts to a map $$H_0(A,A) \To HC_0^-(A).$$ There is an
explicit formula for this lift. For $p \geqslant 3 $ we have
$$ a \mapsto a^p + \sum_{\stackrel{n\ \even,\ p-3\geqslant n \geqslant 2}{\sum i_\alpha = p}} (\mathrm{coefficients}) a^{i_0}
\otimes \cdots \otimes a^{i_n} u^{\frac{n-1}{2}} +
\left( \frac{p-1}{2} \right) !\   a^{\otimes p} u^{\frac{p-1}{2}}, $$
where the last coefficient is non-zero.

The formula for $p=2$ reads: $$a\mapsto a^2+1\otimes a\otimes a \cdot u.$$

\subsection{} One can calculate various simple examples and this seems to suggest a potential mechanism for the degeneration
of the Hodge-to-de Rham spectral sequence in characteristic $p>0.$ The situation is radically different for $p=0$. This mechanism works as follows.

Let us consider polynomials in $u$: $$H^\bullet(C^\redu_\bullet(A,A)[u],\partial+uB).$$ There is no obvious spectral sequence
in this case. What we have is a
quasi-coherent sheaf over $\mathbb{A}^1$ with coordinate $u$.

In characteristic $p=0$ we have that this sheaf vanishes when
$u\neq 0$, namely $(C^\redu_\bullet(A,A)[u,u^{-1}],\partial+uB)$
is acyclic (this is an early observation by Connes). So we have
everything concentrated in an infinite-dimensional stalk at $u=0$,
in the form of an infinite Jordan block, plus some finite Jordan
blocks. It certainly looks like nothing resembling a vector
bundle, it is very singular. The degeneration we seek would mean
that we have no finite Jordan blocks. The situation is
unfortunately quite involved.

In contrast, in finite characteristic $p$, it seems to be the case that the cohomology of the complex $(C^\redu_\bullet(A,A)[u],\partial+uB)$ is
actually a coherent sheaf, and it
will look like a vector bundle if the desired spectral sequence degeneration occurs. The following conjecture explains why we
obtain a vector bundle.

\subsection{Conjecture.} Let $A$ be an flat  dga over $\integer_p$. Let $A_0= A \otimes \integer/p\integer$ over
$\integer/p\integer$. Then $(C^\redu_\bullet(A_0,A_0)[u,u^{-1}],\partial+uB)$ is canonically quasi-isomorphic to
$(C^\redu_\bullet(A_0,A_0)[u,u^{-1}],\partial)$ as $\integer/p[u,u^{-1}]$-modules.

In this conjecture there is no finiteness condition at all.

\subsection{} The reason for this is as follows. The complex $C^\redu_\bullet(A,A)$ admits an obvious increasing filtration
$$\mathrm{Fil}_{\leqslant n} = A \otimes (A/1)^{\otimes \leqslant n-1} + 1 \otimes (A/1)^{\otimes  n}.$$ Let $V:=A/1[1]$, on
$\mathrm{gr}_n(\mathrm{Fil})$ we  can write $\partial + B$ as:
$$\xymatrix{
                   V^{\otimes n}  \ar@/_/[rr]_{1-\sigma} && V^{\otimes n} \ar@/_/[ll]_{1+\sigma+\cdots+\sigma^{n-1}}        }$$
where $\sigma$ is the generator of $\integer/n\integer$ (this
would be acyclic in characteristic $0$). On the other hand
$\partial$ is simply: $$ V^{\otimes n} \stackrel{1-\sigma}{\To}
V^{\otimes n}.$$ For any free $\integer/p\integer$-module the
above complex $\mathrm{gr}_n(\mathrm{Fil})$ with differential
$\partial + B$ is acyclic if $(n,p)=1$. At the same time if $n=kp$
the complex is canonically isomorphic to $$V^{\otimes k}
\stackrel{1-\sigma}{\To} V^{\otimes k}.$$

The hope is that some finite calculation of this sort could allow us to go deep into the spectral sequence and prove the desired degeneration.

\subsection{} We  finish by making some remarks regarding the Weil conjecture. Let $A$ be over $\integer$. It would be
reasonable to hope that we can define
 local $L$-factors by $$L_p(s)=\det \left( 1-\frac{\Fr_p}{p^s}\right)^{- 1},$$ and we should get a sort of non-commutative $L$-function.

We define for a saturated space $X$ an $L$-function: $$L(X)=\prod L_p(s).$$ This $L(X)$ should satisfy
\begin{itemize}
\item The Riemann hypothesis. Namely its zeroes lie on $\Re(s) = \frac{1}{2}$.
\item The Beilinson conjectures. They state that the vanishing order and leading coefficients at $s\in \frac{1}{2}\integer$,
$s\leqslant 1/2 $ are expressed via $K_{1-2s}(X)$ that should in turn be finite dimensional.
\end{itemize}

\subsection{} This $L$-function differs from the traditional $L$-function
defined as a product over all points of a variety over  finite
fields $\mathbb{F}_q$ weighted by $1/q^s$. We can imagine an
$L$-function defined on a saturated space $X$ as the sum over
objects of $\CC_X$ weighted in some way. We do not know the exact
form of these weights but we expect them to depend on certain
stability condition. Such sums appear in string theory as sums
over $D$-branes (for example in the calculation of the entropy of
a black hole \cite{Strominger}). It is reasonable to imagine that
there is a big non-commutative $L$-function one of whose limiting
cases is arithmetic while the other is topological string theory.

\section*{}{\bf Acknowledgements}: I am grateful to Giuseppe Dito, Yan Soibelman and Daniel Sternheimer 
for their help and comments.

\bibliographystyle{amsplain}
\bibliography{Lefschetz}

\end{document}